\documentclass[a4paper,twoside,11pt]{article}
\usepackage{bbm}
\usepackage{amsfonts}
\usepackage{amssymb}\usepackage{amsmath}
 \numberwithin{equation}{section}
\begin{document}

\title{On the Cauchy problem of a two-component b-family equation}

\author{Jingjing Liu\footnote{E-mail:
jingjing830306@163.com,~mcsyzy@mail.sysu.edu.cn },\quad Zhaoyang Yin\\
Department of Mathematics,
Sun Yat-sen University,\\
510275 Guangzhou, China
\bigskip\\
}
\date{}
\maketitle

\begin{abstract}
In this paper, we study the Cauchy problem of a two-component
b-family equation. We first establish the local well-posedness for a
two-component b-family equation by Kato's semigroup theory. Then, we
derive precise blow-up scenarios for strong solutions to the
equation. Moreover, we present several blow-up results for strong
solutions to
the equation.\\


\end{abstract}

\section{Introduction}
\noindent

In this paper we consider the following two-component b-family
equation:
\begin{equation}\label{eq:original}
\left\{\begin{array}{ll}m_{t}=um_{x}+k_{1}u_{x}m+k_{2}\rho\rho_{x},&t > 0,\,x\in \mathbb{R},\\
 \rho_{t}=k_{3}(u\rho )_{x}, &t > 0,\,x\in \mathbb{R},\\
m(0,x) =m_{0}(x),&x\in \mathbb{R}, \\
\rho(0,x) = \rho_{0}(x),&x\in \mathbb{R},\end{array}\right.
\end{equation}
where $m=u-u_{xx}$ and there  are two cases about this system: $(i)$
$k_{1}=b, k_{2}=2b$ and $k_{3}=1$; $(ii)$ $k_{1}=b+1, k_{2}=2$ and
$k_{3}=b$ with $b\in \mathbb{R}$. Eq.(1.1) was recently introduced
by Guha in \cite{p-g}.  The two-component b-family equation is
 defined on a infinite-dimensional Lie group in \cite{M-R}, which is the group
 of orientation-preserving diffeomorphisms of the circle. The group $Diff(\mathbb{S}^{1})$ of smooth
 orientation-preserving diffeomorphisms of the circle $\mathbb{S}^{1}$
is endowed with a smooth manifold structure based on the Fr\'{e}chet
 space $C^{\infty}(\mathbb{S}^{1})$. The composition and inverse are both smooth maps so
that $Diff(\mathbb{S}^{1})$ is a Lie group modeled on Fr\'{e}chet
space, see \cite{p-g} for details.\\

For $\rho\equiv 0$, Eq.(1.1) becomes the b-family equation
\begin{equation}
u_{t}-\alpha^{2}u_{txx}+c_{0}u_{x}+(b+1)uu_{x}+\Gamma
u_{xxx}=\alpha^{2}(bu_{x}u_{xx}+uu_{xxx}).
\end{equation}
Eq.(1.2) can be derived as the family of asymptotically equivalent
shallow water wave equations that emerge at quadratic order accuracy
for any $b\neq -1$ by an appropriate Kodama transformation,
cf.[3-4]. For the case $b=-1$, the corresponding Kodama
transformation is singular and the asymptotic ordering is violated,
cf.[3-4].

With $\alpha=0$ and $b=2$ in Eq.(1.2), we find the well-known KdV
equation which describes the unidirectional propagation of waves at
the free surface of shallow water under the influence of gravity
\cite{D-G-H1}. The Cauchy problem of the KdV equation has been
studied by many authors [6-8] and a satisfactory local or global (in
time) existence theory is now available (e.g. see [7-8]). For $b=2$
and $\gamma=0$, Eq.(1.2) becomes the Camassa-Holm equation,
modelling the unidirectional propagation of shallow water waves over
a flat bottom. The Cauchy problem of the Camassa-Holm equation has
been the subject of a number of studies, for example [9-10]. For
$b=3$  and $c_{0}=\gamma=0$, then we find the Degasperis-Procesi
equation \cite{D-P} from Eq.(1.2), which is regarded as a model for
nonlinear shallow water dynamics. There are also many papers
involving Degasperis-Procesi equation, e.g.[12-13]. The advantage of
the Camassa-Holm equation and the Degasperis-Procesi equation in
comparison with the KdV equation lies in the fact that these two
equations have peakon solitons and models wave breaking [14-15].\\

In \cite{E-Y1} and \cite{Z-Y}, the authors studied Eq.(1.2) on the
line and on the circle respectively for $\alpha> 0$ and $b, c_{0},
\Gamma \in \mathbb{R}$. In \cite{E-Y1} and \cite{Z-Y}, the authors
established the local well-posedness, described the precise blow-up
scenario, proved the equation has strong solutions which exist
globally in time and blow up in finite time. Moreover, the authors
showed the existence of global weak solution to Eq.(1.2) on the line
and on the circle respectively.\\

For $\rho\not\equiv 0$, if $k_{1}=2$, Eq.(1.1) becomes two-component
Camassa-Holm equation. A classical two-component Camassa-Holm
equation
$$
\left\{\begin{array}{ll}m_{t}+um_{x}+2u_{x}m+\sigma\rho\rho_{x}=0,&t > 0,\,x\in \mathbb{R},\\
 \rho_{t}+(u\rho )_{x}=0, &t > 0,\,x\in \mathbb{R},\end{array}\right.
$$
where $m=u-u_{xx}$, $\sigma=\pm 1$ was derived by Constantin and
Ivanov \cite{C-I} in the context of shallow water theory. The
variable $u(x,t)$ describes the horizontal velocity of the fluid and
the variable $\rho(x,t)$ is in connection with the horizontal
deviation of the surface from equilibrium, all measured in
dimensionless units \cite{C-I}. The extended $N=2$ super-symmetric
Camassa-Holm equation was presented recently by Popowicz in
\cite{Po}. The mathematical properties of the two-component
Camassa-Holm equation have been studied in many works cf.[18,
20-25].\\

For $\rho\not\equiv 0$ and $b\in\mathbb{R}$, the Cauchy problem of
Eq.(1.1) has not been studied yet. The aim of this paper is to
establish the local well-posedness, to derive precise blow-up
scenarios, to prove the existence of  strong solutions which blow up
in finite time for Eq.(1.1).\\

Our paper is organized as follows. In Section 2, we establish the
local well posedness of Eq.(1.1). In Section 3, we derive two
precise blow-up scenarios for Eq.(1.1). In Section 4, we discuss the
blow-up phenomena of Eq.(1.1).\\
\newline
\textbf{Notation}  Given a Banach space $Z$, we denote its norm by
 $\|\cdot\|_{Z}$. Since all space of functions are over
 $\mathbb{R}$, for simplicity, we drop $\mathbb{R}$ in our notations of function spaces
  if there is no ambiguity. We let $[A,B]$ denote
 the commutator of linear operator $A$ and $B$. For convenience, we
 let $(\cdot|\cdot)_{s\times r}$ and $(\cdot|\cdot)_{s}$ denote the
 inner products of $H^{s}\times H^{r}$, $s,r\in \mathbb{R}_{+}$ and
 $H^{s}$, $s\in \mathbb{R}_{+}$, respectively.

\section{Local well-posedness}
\noindent

In this section, we will apply Kato's theory to establish the local
well-posedness for the Cauchy problem of Eq.(1.1).\\

\par
For convenience, we state here Kato's theory in the form suitable
for our purpose. Consider the abstract quasi-linear equation:
\begin{equation}
\frac{dv}{dt}+A(v)v=f(v),\ \ t>0,\ \ v(0)=v_{0}.
\end{equation}

\par
Let $X$ and $Y$ be Hilbert spaces such that $Y$ is continuously and
densely embedded in $X$ and let $Q : Y \rightarrow X$ be a
topological isomorphism. Let $L(Y,X)$ denotes the space of all
bounded linear operator from $Y$ to $X$
($L(X)$, if $X=Y$.). Assume that: \\
(i) $A(y)\in L(Y,X)$ for $y\in Y$ with
$$ \|(A(y)-A(z))w\|_{X}\leq\mu_{1}\|
y-z\|_{X}\|w\|_{Y},\ \ \ \ y,z,w\in Y,$$ and $A(y)\in G(X,1,\beta)$,
(i.e. $A(y)$ is quasi-m-accretive),
uniformly on bounded sets in $Y$.\\
(ii) $QA(y)Q^{-1}=A(y)+B(y)$, where $B(y)\in L(X)$ is bounded,
uniformly on bounded sets in $Y$. Moreover,
$$ \|(B(y)-B(z))w\|_{X}\leq\mu_{2}\|
y-z\|_{Y}\|w\|_{X},\ \ \ \ y,z\in Y , \ w\in X.$$ (iii) $f:
Y\rightarrow Y$ and extends also to a map from $X$ to $X$. $f$ is
bounded on bounded sets in $Y$, and
$$ \|f(y)-f(z)\|_{Y}\leq \mu_{3}\|
y-z\|_{Y}, \ \ \ y,z\in Y,$$
$$
\|f(y)-f(z)\|_{X}\leq\mu_{4}\|y-z\|_{X}, \ \ \ y,z\in Y.$$ Here
$\mu_{1}, \mu_{2}, \mu_{3}$ and $ \mu_{4}$ depend only on
max$\{\|y\|_{Y}, \|z\|_{Y} \}$.\\
\newline
\textbf{Theorem 2.1} ~(\cite{K1}). Assume that (i), (ii)
and (iii) hold. Given  $ v_{0}\in Y$,  there is a maximal $T>0$
 depending only on $\parallel v_{0}\parallel_{Y}$ and a
unique solution $v$ to Eq.(2.1) such that
$$
v=v(\cdot,v_{0})\in C([0,T);Y)\cap C^{1}([0,T);X).
$$
Moreover, the map $v_{0}\rightarrow v(\cdot,v_{0})$ is continuous
from Y to
$$
 C([0,T);Y)\cap C^{1}([0,T);X).
$$

\par

We now provide the framework in which we shall reformulate Eq.(1.1).
Note that if $p(x):=\frac{1}{2}e^{-|x|}$, $x\in \mathbb{R}$, then
$(1-
\partial^{2}_{x})^{-1}f = p\ast f $ for all $f \in L^{2}$ and
$p\ \ast \ m=u $. Here we denote by $\ast$ the convolution. Using
this two identities, we can rewrite Eq.(1.1) as follows:
\begin{equation}
\left\{\begin{array}{ll} u_{t}-uu_{x}=\partial_{x}p\ast
(\frac{k_{1}}{2}u^{2}+\frac{3-k_{1}}{2}u_{x}^{2}+\frac{k_{2}}{2}\rho^{2}),&t
> 0,\,x\in \mathbb{R},\\
\rho_{t}-k_{3}u\rho_{x}=k_{3}u_{x}\rho ,&t
> 0,\,x\in \mathbb{R},\\
 u(0,x) = u_{0}(x),&x\in
\mathbb{R},\\
\rho(0,x)=\rho_{0}(x),\ &x\in \mathbb{R}, \end{array}\right.
\end{equation}
or the equivalent form:
\begin{equation}
\left\{\begin{array}{ll}
u_{t}-uu_{x}=\partial_{x}(1-\partial_{x}^{2})^{-1}
(\frac{k_{1}}{2}u^{2}+\frac{3-k_{1}}{2}u_{x}^{2}+\frac{k_{2}}{2}\rho^{2}),&t
> 0,\,x\in \mathbb{R},\\
\rho_{t}-k_{3}u\rho_{x}=k_{3}u_{x}\rho ,&t
> 0,\,x\in \mathbb{R},\\
 u(0,x) = u_{0}(x),&x\in
\mathbb{R},\\
\rho(0,x)=\rho_{0}(x),\ &x\in \mathbb{R}. \end{array}\right.
\end{equation}\\
The main result in this section is the following theorem.\\
\newline
\textbf{Theorem 2.2.}
Given $z_{0}=\left(\begin{array}{c}
                                u_{0} \\
                                \rho_{0} \\
                              \end{array}
                            \right)\in H^{s}\times
H^{s-1}, s\geq 2,$ there exists a maximal $T=T(\parallel
z_{0}\parallel_{H^{s}\times H^{s-1}})>0$, and a unique solution
$z=\left(\begin{array}{c}
                                u \\
                                \rho \\
                              \end{array}
                            \right)$ to
Eq.(2.3) such that
$$
z=z(\cdot,z_{0})\in C([0,T); H^{s}\times H^{s-1})\cap
C^{1}([0,T);H^{s-1}\times H^{s-2}).
$$
Moreover, the solution depends continuously on the initial data,
i.e. the mapping $$z_0\rightarrow z(\cdot,z_0):
 H^{s}\!\times\!
H^{s-1}\!\rightarrow\! C([0,T); H^{s}\!\times\! H^{s-1})\cap
C^{1}([0,T);H^{s-1}\!\times\! H^{s-2})
$$
is continuous.\\

The remainder of this section is devoted to the proof of Theorem
2.2.\\

Let $z:=\left(
          \begin{array}{c}
            u \\
            \rho \\
          \end{array}
        \right)
,\quad A(z)=\left(
         \begin{array}{cc}
           -u\partial_{x} & 0 \\
           0 & -k_{3}u\partial_{x} \\
         \end{array}
       \right)
$  and  $$f(z)=\left(
              \begin{array}{c}
                \partial_{x}(1-\partial_{x}^{2})^{-1}
(\frac{k_{1}}{2}u^{2}+\frac{3-k_{1}}{2}u_{x}^{2}+\frac{k_{2}}{2}\rho^{2}) \\
                k_{3}u_{x}\rho \\
              \end{array}
            \right)
.$$

Set $Y=H^{s}\times H^{s-1},$ $X=H^{s-1}\times H^{s-2}$,
$\Lambda=(1-\partial_{x}^{2})^{\frac{1}{2}}$ and $Q=\left(
                                                     \begin{array}{cc}
                                                       \Lambda & 0 \\
                                                       0 & \Lambda \\
                                                     \end{array}
                                                   \right)
$. Obviously, Q is an isomorphism of $H^{s}\times H^{s-1}$ onto
$H^{s-1}\times H^{s-2}$. In order to prove Theorem 2.2, in view of
Theorem 2.1, we only need to verify $A(z)$ and $f(z)$ satisfy the
conditions (i)-(iii).\\

We first recall the following lemma.\\
\newline
\textbf{Lemma 2.1}~(\cite{K}).
Let r,t be real numbers such that $-r<t\leq r$. Then
$$
\|fg \|_{H^{t}}\leq c\| f\|_{H^{r}}\| g\|_{H^{t}}, \ \ \ \ if \
r>\frac{1}{2},\\
$$
$$
\|fg \|_{H^{t+r-\frac{1}{2}}}\leq c\| f\|_{H^{r}}\|g\|_{H^{t}},  \
if \ r<\frac{1}{2},
$$
where c is a positive  constant depending on r, t.\\

Similar to the proofs of Lemmas 2.5-2.7 in \cite{E-Le-Y}, we get the
following three lemmas.\\
\newline
\textbf{Lemma 2.2.}~
The operator $A(z)=\left(
         \begin{array}{cc}
           -u\partial_{x} & 0 \\
           0 & -k_{3}u\partial_{x} \\
         \end{array}
       \right)
$ with $z\in H^{s}\times H^{s-1},\ s\geq 2$, belongs to
$G(H^{s-1}\times H^{s-2},1,\beta)$.\\
\newline
\textbf{Lemma 2.3.}~
Let $A(z)=\left(
         \begin{array}{cc}
           -u\partial_{x} & 0 \\
           0 & -k_{3}u\partial_{x} \\
         \end{array}
       \right)
$  with $z\in H^{s}\times H^{s-1},\ s\geq 2$. Then $A(z)\in
L(H^{s}\times H^{s-1},H^{s-1}\times H^{s-2})$ and
$$
\|(A(z)-A(y))w\|_{H^{s}\times H^{s-1}}\leq \mu_{1}\|
z-y\|_{H^{s-1}\times H^{s-2}}\| w\|_{H^{s}\times H^{s-1}}, $$ for
all $z, y, w\in H^{s}\times H^{s-1}$.\\
\newline
\textbf{Lemma 2.4.}~
Let $B(z)=QA(z)Q^{-1}-A(z)$ with $z\in H^{s}\times H^{s-1},$
$s\geq2.$ Then $B(z)\in L(H^{s-1}\times H^{s-2})$ and
$$
\|(B(z)-B(y))w\|_{H^{s-1}\times H^{s-2}}\leq
\mu_{2}\|z-y\|_{H^{s}\times H^{s-1}}\parallel w\|_{H^{s-1}\times
H^{s-2}}, $$ for all $z, y\in H^{s}\times H^{s-1}$ and $w\in
H^{s-1}\times H^{s-2}$.\\

We now prove that $f$ satisfies the condition (iii) in Theorem 2.1.\\
\newline
\textbf{Lemma 2.5.}~
Let $z\in H^{s}\times H^{s-1}, s\geq 2$ and let $$f(z)=\left(
              \begin{array}{c}
                \partial_{x}(1-\partial_{x}^{2})^{-1}
(\frac{k_{1}}{2}u^{2}+\frac{3-k_{1}}{2}u_{x}^{2}+\frac{k_{2}}{2}\rho^{2}) \\
                k_{3}u_{x}\rho \\
              \end{array}
            \right)
.$$ Then f is bounded on bounded sets in $H^{s}\times H^{s-1}$, and
for all $y,z\in H^{s}\times H^{s-1}$ satisfies
$$
(a)\quad \|f(y)-f(z)\|_{H^{s}\times H^{s-1}}\leq\mu_{3}\|
y-z\|_{H^{s}\times H^{s-1}},$$
$$
(b)\quad \|f(y)-f(z)\|_{H^{s-1}\times H^{s-2}}\leq\mu_{3}\|
y-z\|_{H^{s-1}\times H^{s-2}} .$$
\newline
\textbf{Proof}~
Let $y,z\in {H^{s}\times H^{s-1}}, s\geq 2$. Note that $H^{s-1}$ is
a Banach algebra. Then, we have
\begin{eqnarray*}
&& \|f(y)-f(z)\|_{H^{s}\times
H^{s-1}}\\
&\leq & \|\partial_{x}(1-\partial_{x}^{2})^{-1}
(\frac{k_{1}}{2}(y_{1}^{2}-u^{2})-\frac{3-k_{1}}{2}(y_{1,x}^{2}-u_{x}^{2})+\frac{k_{2}}{2}(y_{2}^{2}-\rho^{2})\|_{H^{s}}
+|k_{3}|\\&&
\|y_{1,x}y_{2}-u_{x}\rho\|_{H^{s-1}}\\
&\leq&\frac{|k_{1}|}{2}\|(y_{1}-u)(y_{1}+u)\|_{H^{s-1}}+\frac{|3-k_{1}|}{2}\|(y_{1,x}-u_{x})(y_{1,x}+u_{x})\|_{H^{s-1}}
+\\&&
\frac{|k_{2}|}{2}\|(y_{2}-\rho)(y_{2}+\rho)\|_{H^{s-1}}+|k_{3}|(\|y_{1,x}y_{2}-y_{1,x}\rho\|_{H^{s-1}}+\|y_{1,x}\rho-u_{x}\rho\|_{H^{s-1}})\\
&\leq&\frac{|k_{1}|}{2}\|y_{1}-u\|_{H^{s-1}}\|y_{1}+u\|_{H^{s-1}}+\frac{|3-k_{1}|}{2}\|y_{1}-u\|_{H^{s}}\|y_{1}+u\|_{H^{s}}\\&&
+\frac{|k_{1}|}{2}\|y_{2}-\rho\|_{H^{s-1}}\|y_{2}+\rho\|_{H^{s-1}}+|k_{3}|\|y_{1}\|_{H^{s}}\|y_{2}-\rho\|_{H^{s-1}}
\\&&
+|k_{3}|\|y_{1}-u\|_{H^{s}}\|\rho\|_{H^{s-1}}\\
&\leq&
\left(\frac{|k_{1}|}{2}+\frac{|3-k_{1}|}{2}+\frac{|k_{2}|}{2}+2|k_{3}|\right)(\|y\|_{H^{s}\times
H^{s-1}}+\|z\|_{H^{s}\times H^{s-1}})\|y-z\|_{H^{s}\times H^{s-1}}.
\end{eqnarray*}
 This proves (a). Taking $y=0$ in the above inequality, we obtain that $f$ is
bounded on bounded set in $H^s\times H^{s-1}$.\\
\par
Next, we prove
(b). Note that $H^{s-1}$ is a Banach algebra. Then, we have
\begin{eqnarray*}
&& \|f(y)-f(z)\|_{H^{s-1}\times
H^{s-2}}\\
&\leq & \|\partial_{x}(1-\partial_{x}^{2})^{-1}
(\frac{k_{1}}{2}(y_{1}^{2}-u^{2})-\frac{3-k_{1}}{2}(y_{1,x}^{2}-u_{x}^{2})+\frac{k_{2}}{2}(y_{2}^{2}-\rho^{2})\|_{H^{s-1}}
+\\&&
|k_{3}|\|y_{1,x}y_{2}-u_{x}\rho\|_{H^{s-2}}\\
&\leq&\frac{|k_{1}|}{2}\|(y_{1}-u)(y_{1}+u)\|_{H^{s-2}}+\frac{|3-k_{1}|}{2}\|(y_{1,x}-u_{x})(y_{1,x}+u_{x})\|_{H^{s-2}}+\\&&
\frac{|k_{2}|}{2}\|(y_{2}-\rho)(y_{2}+\rho)\|_{H^{s-2}}+|k_{3}|\|y_{1,x}y_{2}-y_{1,x}\rho\|_{H^{s-2}}+|k_{3}|\|y_{1,x}\rho-u_{x}\rho\|_{H^{s-2}}\\
&\leq&\frac{|k_{1}|}{2}\|y_{1}-u\|_{H^{s-2}}\|y_{1}+u\|_{H^{s-2}}+\frac{c|3-k_{1}|}{2}\|y_{1,x}-u_{x}\|_{H^{s-2}}\|y_{1,x}
+u_{x}\|_{H^{s-1}}\\&&
+\frac{c|k_{2}|}{2}\|y_{2}-\rho\|_{H^{s-2}}\|y_{2}+\rho\|_{H^{s-1}}+c|k_{3}|\|y_{1,x}\|_{H^{s-1}}\|y_{2}-\rho\|_{H^{s-2}}
\\&&+c|k_{3}|\|y_{1}-u\|_{H^{s-1}}\|\rho\|_{H^{s-1}}\\
&\leq&
\left(\frac{|k_{1}|}{2}+\frac{c|3-k_{1}|}{2}+\frac{c|k_{2}|}{2}+2c|k_{3}|\right)(\|y\|_{H^{s}\times
H^{s-1}}+\|z\|_{H^{s}\times H^{s-1}})\|y-z\|_{H^{s-1}\times
H^{s-2}},
\end{eqnarray*}
where we applied Lemma 2.1 with $r=s-1$, $t=s-2$ and $c$ only
depends on $k_{1}, k_{2}$ and $k_{3}$. This proves (b) and completes
the proof of the lemma.\\

\noindent\textit{\ Proof of Theorem 2.2:} Combining Theorem 2.1 and
Lemmas 2.2-2.5, we get the statement of Theorem 2.2.

\section{Precise blow-up scenarios}
\noindent

In this section, we will derive precise blow-up scenarios for strong
solutions to Eq.(1.1).\\

We first recall the following two useful lemmas.\\
\newline
\textbf{Lemma 3.1}~
(\cite{K1}). If $r>0$, then $H^{r}\cap L^{\infty}$ is an algebra.
Moreover
$$\|fg\|_{H^r}\leq
c(\|f\|_{L^{\infty}}\|g\|_{H^r}+\|f\|_{H^r}\|g\|_{L^{\infty}}),
$$
where c is a constant depending only on r.\\
\newline
\textbf{Lemma 3.2}~
(\cite{K1}). If $r>0$, then
$$\|[\Lambda^{r},f]g\|_{L^{2}}\leq c(\|\partial_{x}f\|_{L^{\infty}}
\|\Lambda^{r-1}g\|_{L^2}+\|\Lambda^r f\|_{L^2}\| g\|_{L^{\infty}}),
$$
where c is a constant depending only on r.\\

Then, we have the following useful result.\\
\newline
\textbf{Theorem 3.1.}~
 Let $z_{0}=\left(\begin{array}{c}
                                u_{0} \\
                                \rho_{0} \\
                              \end{array}
                            \right)
\in H^{s}\times H^{s-1}$, $s\geq2,$ be given and assume that T is
the maximal existence time of the corresponding solution
                  $z=\left(\begin{array}{c}
                                      u \\
                                      \rho \\
                                    \end{array}
                                  \right)
$ to Eq.(1.1) with the initial data $z_{0}$. If there exists $M>0$
such that
$$\|u_{x}(t,\cdot)\|_{L^{\infty}}+\|\rho(t,\cdot)\|_{L^{\infty}}
+\|\rho_{x}(t,\cdot)\|_{_{L^{\infty}}}\leq M,\ \ t\in[0,T),
$$
then the $H^s\times H^{s-1}$-norm of $z(t,\cdot)$ does not blow up
on [0,T).\\
\newline
\textbf{Proof}~
Let $z=\left(
                                    \begin{array}{c}
                                      u \\
                                      \rho \\
                                    \end{array}
                                  \right)$ be the solution to
Eq.(1.1) with the initial data $z_{0}\in H^s\times H^{s-1},\ s\geq
2$, and let T be the maximal existence time of the corresponding
solution $z$, which is guaranteed by Theorem 2.2.\\

Applying the operator $\Lambda^{s}$ to the first equation in (2.3),
multiplying by $\Lambda^{s} u$, and integrating over $\mathbb{R}$,
we obtain
\begin{eqnarray}
&&\frac{d}{dt}\|u\|^{2}_{H^{s}}\\&=&\nonumber2(uu_{x},
u)_{s}+2(u,f_{11}(u))_{s}+2(u,f_{12}(u))_{s},
\end{eqnarray}
where
$$f_{11}(u)=\partial_{x}(1-\partial_{x}^{2})^{-1}\frac{k_{1}}{2}u^{2}=(1-\partial_{x}^{2})^{-1}(k_{1}uu_{x}),$$
and
$$f_{12}(u)=\partial_{x}(1-\partial_{x}^{2})^{-1}(\frac{3-k_{1}}{2}u_{x}^{2}+\frac{k_{2}}{2}\rho^{2}).$$From
the proof of Theorem 3.1 in \cite{E-Le-Y}, we have
$$
|(uu_{x},u)_{s}| \leq c\| u_{x}\|_{L^{\infty}}\|
 u\|^2_{H^s},$$
and
$$|(f_{11}(u),u)_{s}|\leq
c\|u_{x}\|_{L^{\infty}}\|u\|_{H^{s}}^{2}.$$ Furthermore, we estimate
the third term of the right hand side of Eq.(3.1) in the following
way:
\begin{eqnarray*}
&&|(f_{12}(u),u)_{s}|\\
&\leq&\|f_{12}(u)\|_{H^{s}}\|u\|_{H^{s}}\\
&=&\|\partial_{x}(1-\partial_{x}^{2})^{-1}(\frac{3-k_{2}}{2}u_{x}^{2}+\frac{k_{2}}{2}\rho^{2})\|_{H^{s}}\|u\|_{H^{s}}\\
&\leq&\left(\frac{|3-k_{1}|}{2}\|u_{x}^{2}\|_{H^{s-1}}+\frac{|k_{2}|}{2}\|\rho^{2}\|_{H^{s-1}}\right)\|u\|_{H^{s}}\\
&\leq&(|3-k_{1}|c\|u_{x}\|_{L^{\infty}}\|u_{x}\|_{H^{s-1}}+|k_{2}|c\|\rho\|_{L^{\infty}}\|\rho\|_{H^{s-1}})\|u\|_{H^{s}}\\
&\leq&c(\|u_{x}\|_{L^{\infty}}+\|\rho\|_{L^{\infty}})(\|\rho\|_{H^{s-1}}^{2}+\|u\|_{H^{s}}^{2}).
\end{eqnarray*}
Here, we applied Lemma 3.1 with $r=s-1$. Combining the above three
inequalities with (3.1), we get
\begin{equation}
\frac{d}{dt}\|u\|^2_{H^s}\leq c(\|u_{x}\|_{L^{\infty}}
 +\|\rho\|_{L^{\infty}})(\|\rho\|_{H^{s-1}}^2+\|
u\|_{H^s}^2).
\end{equation}
In order to derive a similar estimate for the second component
$\rho$, we apply the operator $\Lambda^{s-1}$ to the second equation
in (2.3), multiply by $\Lambda^{s-1}\rho$, and integrate over
$\mathbb{R}$ we obtain
\begin{equation*}
\frac{d}{dt}\|
\rho\|^2_{H^{s-1}}=2k_{3}(u\rho_{x},\rho)_{s-1}+2k_{3}(u_{x}\rho,\rho)_{s-1}.
\end{equation*}
Following the similar argument in the proof of Theorem 3.1 in
\cite{E-Le-Y}, we have
$$ |(u\rho_{x},\rho)_{s-1}|\leq c(\|u_{x}\|_{L^{\infty}}
 +\|\rho_{x}\|_{L^{\infty}})(\|\rho\|_{H^{s-1}}^2+\|u\|_{H^s}^2)
$$
and
$$ |(u_{x}\rho,\rho)_{s-1}|\leq c(\|u_{x}\|_{L^{\infty}}
 +\|\rho\|_{L^{\infty}})(\|\rho\|_{H^{s-1}}^2+\|u\|_{H^s}^2).
$$
Then, it follows that
\begin{eqnarray}
&&\frac{d}{dt}\|\rho\|^2_{H^{s-1}}\\&\leq&c(\|u_{x}\|_{L^{\infty}}+\|\rho_{x}\|_{L^{\infty}}
 +\|\rho\|_{L^{\infty}})(\|\rho\|_{H^{s-1}}^2+\|u\|_{H^s}^2)\nonumber.
\end{eqnarray}
By (3.2) and (3.3), we obtain
\begin{eqnarray*}
&&\frac{d}{dt}(\|u\|^2_{H^{s}}+\|\rho\|^2_{H^{s-1}})\\&\leq & c(\|
u_{x}\|_{L^{\infty}}+\|\rho_{x}\|_{L^{\infty}}
 +\|\rho\|_{L^{\infty}})(\|\rho\|_{H^{s-1}}^2+\|u\|_{H^s}^2).
\end{eqnarray*}
An application of Gronwall's inequality and the assumption of the
theorem yield
$$
\| u\|^2_{H^{s}}+\|\rho\|^2_{H^{s-1}}\leq \exp(cMt)(\|
u_{0}\|^2_{H^{s}}+\|\rho_{0}\|^2_{H^{s-1}}).
$$
This completes the proof of the theorem.\\

Consider now the following initial value problem
\begin{equation}
\left\{\begin{array}{ll}q_{t}=u(t,-k_{3}q),\ \ \ \ t\in[0,T), \\
q(0,x)=x,\ \ \ \ x\in\mathbb{R}, \end{array}\right.
\end{equation}
where $u$ denotes the first component of the solution $z$ to
Eq.(1.1). Applying classical results in the theory of ordinary
differential equations, one can obtain two results on $q$ which are
crucial in studying  blow-up phenomena.\\
\newline
\textbf{Lemma 3.3.}~Let $u\in C([0,T);H^s)\cap
C^1([0,T);H^{s-1}), s\geq 2$. Then Eq.(3.4) has a unique solution
$q\in C^1([0,T)\times \mathbb{R};\mathbb{R})$. Moreover, the map
$q(t,\cdot)$ is an increasing diffeomorphism of $\mathbb{R}$ with
$$
q_{x}(t,x)=\exp\left(\int_{0}^{t}-k_{3}u_{x}(s,-k_{3}q(s,x))ds\right)>
0, \ \ \forall(t,x)\in [0,T)\times \mathbb{R}.$$
\newline
\textbf{Proof}
Since $u\in C([0,T);H^s)\cap C^1([0,T);H^{s-1}), s\geq 2$ and
$H^{s}\subset C^{1},$ we see that both functions $u(t,x)$ and
$u_{x}(t,x)$ are bounded, Lipschitz in the space variable $x$, and
of class $C^{1}$ in time. Therefore, for fixed $x\in \mathbb{R}$,
equation (3.4) is an ordinary differential equation. Then well-known
classical results in the theory of ordinary differential equations
yield that equation (3.4) has a unique solution $q\in
C^1([0,T)\times \mathbb{R};\mathbb{R})$.

Differentiation of equation (3.4) with respect to $x$ yields
\begin{equation}
\left\{\begin{array}{ll}\frac{d}{dt}q_{x}=-k_{3}u_{x}(t,-k_{3}q)q_{x},\ \ \ \ t\in[0,T), \\
q(0,x)=x,\ \ \ \ \ \ \ \ \ \  \ \ \ \ \ \ \ x\in\mathbb{R}.
\end{array}\right.
\end{equation}
The solution to equation (3.5) is given by
\begin{equation}
q_{x}(t,x)=\exp\left(\int_{0}^{t}-k_{3}u_{x}(s,-k_{3}q(s,x))ds\right),
\ \ \forall(t,x)\in [0,T)\times \mathbb{R}.
\end{equation} For every
$T^{\prime}< T$, by Sobolev¡¯s imbedding theorem, we
get$$\sup\limits_{(s,x)\in [0,T^{\prime})\times
\mathbb{R}}|u_{x}(s,x)|< \infty.$$ Thus, we infer from Eq.(3.6) that
there exists a constant $K>0$ such that $q_{x}(t,x)\geq
e^{-|k_{3}|Kt}>0$ for $(t,x)\in [0,T^{\prime})\times \mathbb{R}.$
This completes the proof of the lemma.\\
\newline
\textbf{Lemma 3.4.}~ Let $z_{0}\in H^s\times H^{s-1}, s\geq 2$, and let $T>0$ be
the maximal existence time of the corresponding solution
$z=\left(\begin{array}{c}
                                      u \\
                                      \rho \\
                                    \end{array}
                                  \right)
$ to Eq.(2.3). Then we have
\begin{equation}
\rho(t,-k_{3}q(t,x))q_{x}(t,x)=\rho_{0}(-k_{3}x),\ \ \ \forall \
(t,x)\in [0,T)\times \mathbb{R}.
\end{equation}
Moreover, if $k_{3}\leq 0$ and there exists $M>0 $ such that
$u_{x}(t,x)\geq-M$ for all $(t,x)\in [0,T)\times \mathbb{R}$, then
$$
\|\rho(t,\cdot)\|_{L^{\infty}}=\|\rho(t,-k_{3}q(t,\cdot))\|_{L^{\infty}}\leq
e^{-k_{3}MT}\|\rho_{0}(\cdot)\|_{L^{\infty}},\ \ \ \forall \ t\in
[0,T);$$ if  $k_{3}\geq 0$ and there exists $M>0 $ such that
$u_{x}(t,x)\leq M$ for all $(t,x)\in [0,T)\times \mathbb{R}$, then
$$
\|\rho(t,\cdot)\|_{L^{\infty}}=\|\rho(t,-k_{3}q(t,\cdot))\|_{L^{\infty}}\leq
e^{k_{3}MT}\|\rho_{0}(\cdot)\|_{L^{\infty}},\ \ \ \forall \ t\in
[0,T).$$ Furthermore, if there exists $M>0 $ such that
$|u_{x}(t,x)|\leq M$ for all $(t,x)\in [0,T)\times \mathbb{R}$, then
$$
\|\rho(t,\cdot)\|_{L^{\infty}}=\|\rho(t,-k_{3}q(t,\cdot))\|_{L^{\infty}}\leq
e^{|k_{3}|MT}\|\rho_{0}(\cdot)\|_{L^{\infty}},\ \ \ \forall \ t\in
[0,T).$$
\newline
\textbf{Proof}\quad Differentiating the left-hand side of Eq.(3.7) with
respect to $t$, in view of (3.4) and Eq.(2.3), we obtain
\begin{eqnarray*}
&&\frac{d}{dt}{(\rho(t,-k_{3}q(t,x))q_{x}(t,x))}\\&=&(\rho_{t}+\rho_{x}\cdot(-k_{3}q_{t}(t,x)))q_{x}(t,x)
+\rho \cdot q_{xt}(t,x)\\
&=&(\rho_{t}-k_{3}\rho_{x}u)q_{x}(t,x)-k_{3}\rho u_{x} q_{x}(t,x)\\
&=&(\rho_{t}-k_{3}\rho_{x}u-k_{3}\rho u_{x})q_{x}(t,x)=0
\end{eqnarray*}%
This proves (3.7). By Lemma 3.3, in view of (3.7), we obtain
\begin{eqnarray*}
\|\rho(t,\cdot)\|_{L^{\infty}}&=&\|\rho(t,-k_{3}q(t,\cdot))\|_{L^{\infty}}\\
&=&\|\exp\left(k_{3}\int_0^tu_x(s,-k_{3}q(s,\cdot))ds\right)\rho_0(-k_{3}x)\|_{L^{\infty}}.\
\ \forall \ t\in[0,T).
\end{eqnarray*}
The left proof is obvious, so we omit it here.\\

By Theorem 3.1 and Lemma 3.4 we have the following corollary.\\
\newline
\textbf{Corollary 3.1.}
 Let $z_{0}=\left(\begin{array}{c}
                                u_{0} \\
                                \rho_{0} \\
                              \end{array}
                            \right)
\in H^{s}\times H^{s-1}$, $s\geq2$, be given and assume that T is
the maximal existence time of the corresponding solution
                  $z=\left(\begin{array}{c}
                                      u \\
                                      \rho \\
                                    \end{array}
                                  \right)
$ to Eq.(1.1) with the initial data $z_{0}$. If there exists $M>0$
such that
$$\|u_{x}(t,\cdot)\|_{L^{\infty}}
+\|\rho_{x}(t,\cdot)\|_{L^{\infty}}\leq M,\ \ t\in[0,T),
$$
then the $H^s\times H^{s-1}$-norm of $z(t,\cdot)$ does not blow up
on [0,T).\\

Our next result describes the precise blow-up scenario for
sufficiently regular solutions to Eq.(1.1).\\
\newline
\textbf{Theorem 3.2.} ~Let $z_{0}=\left(
                                                     \begin{array}{c}
                                                       u_{0} \\
                                                       \rho_{0} \\
                                                     \end{array}
                                                   \right)
\in H^s\times H^{s-1}, s>\frac{5}{2},$ and let T be the maximal
existence time of the corresponding solution $z=\left(
                                    \begin{array}{c}
                                      u \\
                                      \rho \\
                                    \end{array}
                                  \right)
$ to Eq.(1.1). If $k_{1}\leq -\frac{1}{2},\ k_{3}\leq \min\{0,
k_{2}\},$ then $z$ blows up in finite time if and only if
$$
\liminf_{t\rightarrow T}\inf_{x\in\mathbb{R}}\{u_{x}(t,x)\}=-\infty
\ \ \text{ or }\ \ \limsup_{t\rightarrow
T}\|\rho_{x}(t,\cdot)\|_{L^{\infty}}=+\infty.$$ If $k_{1}\geq 1,\
k_{3}\geq \max\{0, k_{2}\},$ then $z$ blows up in finite time if and
only if
$$
\limsup_{t\rightarrow T}\sup_{x\in\mathbb{R}}\{u_{x}(t,x)\}=+\infty
\ \ \text{ or }\ \ \limsup_{t\rightarrow
T}\|\rho_{x}(t,\cdot)\|_{L^{\infty}}=+\infty.$$ Otherwise, $z$ blows
up in finite time if and only if
$$
\limsup_{t\rightarrow T}\|u_{x}(t,\cdot)\|_{L^{\infty}}=+\infty \ \
\text{ or }\ \ \limsup_{t\rightarrow
T}\|\rho_{x}(t,\cdot)\|_{L^{\infty}}=+\infty.$$
\newline
\textbf{Proof}\quad
Let $z=\left(
                                     \begin{array}{c}
                                       u \\
                                       \rho \\
                                     \end{array}
                                   \right)
$ be the solution to Eq.(1.1) for the initial $z_{0}\in H^s\times
H^{s-1}, s>\frac{5}{2}$, and let $T$ be the maximal existence time
of the solution $z$, which is guaranteed by Theorem 2.2.\\

Multiplying the first equation in (1.1) by $m=u-u_{xx}$ and
integrating by parts, we get
\begin{eqnarray}
&&\frac{d}{dt}\int_{\mathbb{R}}m^2dx\\
\nonumber&=&(2k_{1}-1)\int_{\mathbb{R}}m^{2}u_{x}dx+2k_{2}\int_{\mathbb{R}}u\rho\rho_{x}-2k_{2}\int_{\mathbb{R}}u_{xx}\rho\rho_{xx}\\
\nonumber&=&(2k_{1}-1)\int_{\mathbb{R}}m^{2}u_{x}dx-k_{2}\int_{\mathbb{R}}u_{x}\rho^{2}+k_{2}\int_{\mathbb{R}}u_{xxx}\rho^{2}.
\end{eqnarray}\\

Differentiating the first equation in (1.1) with respect to $x$,
multiplying the obtained equation by $m_{x}=u_{x}-u_{xxx}$, and
integrating by parts, we get
\begin{eqnarray}
&&\frac{d}{dt}\int_{\mathbb{R}}m_{x}^{2}dx\\
\nonumber&=&2(1+k_{1})\int_{\mathbb{R}}u_{x}m_{x}^{2}dx+2\int_{\mathbb{R}}um_{x}m_{xx}+2k_{1}\int_{\mathbb{R}}u_{xx}mm_{x}\\
\nonumber&&+2k_{2}\int_{\mathbb{R}}m_{x}\rho_{x}^{2}+2k_{2}\int_{\mathbb{R}}m_{x}\rho\rho_{xx}\\
\nonumber&=&(1+2k_{1})\int_{\mathbb{R}}u_{x}m_{x}^{2}-k_{1}\int_{\mathbb{R}}u_{x}m^{2}
+k_{2}\int_{\mathbb{R}}u_{xxx}(\rho^{2}-2\rho\rho_{xx}-2\rho_{x}^{2}).
\end{eqnarray}
Here we used the relation $\int_{\mathbb{R}}m^{2}m_{x}dx=0.$\\

Combining (3.8) with (3.9) and integrating by parts, we obtain
\begin{eqnarray}
&&\frac{d}{dt}\int_{\mathbb{R}}m^2+m^2_{x}dx\\
\nonumber&=&(k_{1}-1)\int_{\mathbb{R}}u_{x}m^{2}+(1+2k_{1})\int_{\mathbb{R}}u_{x}m_{x}^{2}
-k_{2}\int_{\mathbb{R}}u_{x}\rho^{2}\\
\nonumber&&+2k_{2}\int_{\mathbb{R}}u_{xxx}(\rho^{2}-\rho\rho_{xx}-\rho_{x}^{2})dx.
\end{eqnarray}

Multiplying the second equation in (1.1) by $\rho$ and integrating
by parts, we deduce
\begin{equation}
\frac{d}{dt}\int_{\mathbb{R}}\rho^2dx=k_{3}\int_{\mathbb{R}}u_{x}\rho^2dx.
\end{equation}\\

Differentiating the second equation in (1.1) with respect to $x$,
multiplying the obtained equation by $\rho_{x}$, and integrating by
parts, we find
\begin{eqnarray}
&&\frac{d}{dt}\int_{\mathbb{R}}\rho^2_{x}dx\\
\nonumber&=&4k_{3}\int_{\mathbb{R}}u_{x}\rho_{x}^{2}+2k_{3}\int_{\mathbb{R}}u\rho_{x}\rho_{xx}+
2k_{3}\int_{\mathbb{R}}u_{xx}\rho\rho_{x}\\
\nonumber&=&3k_{3}\int_{\mathbb{R}}u_{x}\rho_{x}^{2}-k_{3}\int_{\mathbb{R}}u_{xxx}\rho^{2}.
\end{eqnarray}\\

Differentiating the second equation in (1.1) with respect to $x$
twice, multiplying the obtained equation by $\rho_{xx}$, and
integrating by parts, we obtain
\begin{eqnarray}
&&\frac{d}{dt}\int_{\mathbb{R}}\rho^2_{xx}dx\\
\nonumber&=&5k_{3}\int_{\mathbb{R}}u_{x}\rho_{xx}^{2}dx+k_{3}\int_{\mathbb{R}}u_{xxx}(2\rho\rho_{xx}-3\rho_{x}^{2})dx.
\end{eqnarray}

Thus, in view of (3.10)-(3.13), we obtain
\begin{eqnarray}
&&
\frac{d}{dt}\int_{\mathbb{R}}(m^2+m^2_{x}+\rho^2+\rho^2_{x}+\rho^2_{xx})dx\\
\nonumber&=&(k_{1}-1)\int_{\mathbb{R}}u_{x}m^{2}dx+(1+2k_{1})\int_{\mathbb{R}}u_{x}m_{x}^{2}dx
\\
\nonumber&&+(k_{3}-k_{2})\int_{\mathbb{R}}u_{x}\rho^{2}dx+3k_{3}\int_{\mathbb{R}}u_{x}\rho_{x}^{2}dx+5k_{3}\int_{\mathbb{R}}u_{x}\rho_{xx}^{2}dx
\\
\nonumber&&+\int_{\mathbb{R}}u_{xxx}((2k_{2}-k_{3})\rho^{2}+2(k_{3}-k_{2})\rho\rho_{xx}-(2k_{2}+3k_{3})\rho_{x}^{2}dx.
\end{eqnarray}

If $k_{1}\leq -\frac{1}{2},\ k_{3}\leq \min\{0, k_{2}\},$ assume
that there exist $M_{1}>0$ and $M_{2}>0$ such that $$u_{x}(t,x)\geq
- M_{1},\ \ \ \forall \ (t,x)\in [0,T)\times\mathbb{R},$$ and
$$
\|\rho_{x}(t,\cdot)\|_{L^{\infty}}\leq M_{2}, \ \ \ \forall \ t\in
[0,T).$$ By Lemma 3.4, we deduce
$$\|\rho(t,\cdot)\|_{L^{\infty}}\leq
e^{-k_{3}M_{1}T}\|\rho_{0}(\cdot)\|_{L^{\infty}}, \ \ \ \forall \
t\in [0,T).$$ It then follows from (3.14) that
\begin{eqnarray}
&&\frac{d}{dt}\int_{\mathbb{R}}(m^2+m^2_{x}+\rho^2+\rho^2_{x}+\rho^2_{xx})dx\\
&\leq&
(-3k_{1}+k_{2}-9k_{3})M_{1}\int_{\mathbb{R}}(m^2+m^2_{x}+\rho^{2}+\rho^2_{x}+\rho^2_{xx})dx\nonumber
\\&&+(|2k_{2}-k_{3}|\|\rho\|_{L^{\infty}(\mathbb{R})}+2|k_{3}-k_{2}|\|\rho\|_{L^{\infty}(\mathbb{R})}+|2k_{2}+3k_{3}|M_{2})\nonumber\\
&& \int_{\mathbb{R}}(3u^2_{xxx}+
\rho^2+\rho^2_{x}+\rho_{xx}^2)dx\nonumber\\&\leq&c\int_{\mathbb{R}}(m^2+m^2_{x}+\rho^2+\rho^2_{x}+\rho^2_{xx})dx\nonumber,
\end{eqnarray}
where
$$c=(-3k_{1}+k_{2}-9k_{3})M_{1}+3((|2k_{2}-k_{3}|+
2|k_{3}-k_{2}|)e^{-k_{3}M_{1}T}\|\rho_{0}(\cdot)\|_{L^{\infty}}+|2k_{2}+3k_{3}|M_{2}).$$
By means of Gronwall's inequality, we obtain
\begin{eqnarray*}
\|u(t,\cdot)\|^2_{H^3}+\|\rho(t,\cdot)\|^2_{H^2}&\leq&\|m(t,\cdot)\|^2_{H^1}+\|\rho(t,\cdot)\|^2_{H^2}\\
&\leq&(\|m(0,\cdot)\|^2_{H^1}+\|\rho(0,\cdot)\|^2_{H^2 })e^{ct}, \ \
\ \ \ \ \forall \ t\in [0,T).
\end{eqnarray*}
The above inequality, Sobolev's imbedding theorem and Corollary 3.1
ensure that the solution $z$ does not blow up in finite time.\\

If $k_{1}\geq 1,\ k_{3}\geq \max\{0, k_{2}\},$ assume that there
exist $M_{1}>0$ and $M_{2}>0$ such that $$u_{x}(t,x)\leq M_{1},\ \ \
\forall \ (t,x)\in [0,T)\times\mathbb{R},$$ and
$$
\|\rho_{x}(t,\cdot)\|_{L^{\infty}}\leq M_{2}, \ \ \ \forall \ t\in
[0,T).$$ By Lemma 3.4, we deduce
$$\|\rho(t,\cdot)\|_{L^{\infty}}\leq
e^{k_{3}M_{1}T}\|\rho_{0}(\cdot)\|_{L^{\infty}}, \ \ \ \forall \
t\in [0,T).$$ The left proof is similar to the proof above, so we
omit it here.\\

 Otherwise, let $T< \infty.$ Assume that there exists $M_{1}>0$ and $M_{2}>0$ such
that
$$
\|u_{x}(t,\cdot)\|_{L^{\infty}}\leq M_{1},\ \ \ \forall \ t\in
[0,T),$$ and
$$
\|\rho_{x}(t,\cdot)\|_{L^{\infty}}\leq M_{2}, \ \ \ \forall \ t\in
[0,T).$$ We can get
$$\|u_{x}(t,\cdot)\|_{L^{\infty}}
+\|\rho_{x}(t,\cdot)\|_{_{L^{\infty}}}\leq M,\ \ t\in[0,T).
$$ Then Corollary 3.1 implies that the solution
$z$ does not blow up in finite time.\\

On the other hand, by  Sobolev's imbedding theorem, we see that if
one of the conditions in the theorem holds, then the solution will
blow up in finite time. This completes the proof of the theorem.\\

For initial data $z_{0}\in H^{2}\times H^{1}$, we have the following
precise blow-up scenario.\\
\newline
\textbf{Theorem 3.3.}~ Let $z_{0}=\left(
                                                     \begin{array}{c}
                                                       u_{0} \\
                                                       \rho_{0} \\
                                                     \end{array}
                                                   \right)
\in H^2\times H^{1},$ and let T be the maximal existence time of the
corresponding solution $z=\left(
                                    \begin{array}{c}
                                      u \\
                                      \rho \\
                                    \end{array}
                                  \right)
$ to Eq.(1.1). If $k_{1}\leq\frac{1}{2},$ $k_{3}\leq\min\{
k_{2},0\}$, then the corresponding solution blows up in finite time
if and only if
$$
\liminf_{t\rightarrow T}\inf_{x\in\mathbb{R}}\{u_{x}(t,x)\}=-\infty
.$$ If $k_{1}\geq \frac{1}{2},$ $k_{3}\geq\max\{k_{2},0\}$, then the
corresponding solution blows up in finite time if and only if
$$
\limsup_{t\rightarrow T}\sup_{x\in\mathbb{R}}\{u_{x}(t,x)\}=+\infty
.$$ Otherwise, the corresponding solution blows up in finite time if
and only if
$$
\limsup_{t\rightarrow T}\|u_{x}(t,\cdot)\|_{L^{\infty}}=+\infty .$$
\newline
\textbf{Proof}\quad Let $z=\left(
                                     \begin{array}{c}
                                       u \\
                                       \rho \\
                                     \end{array}
                                   \right)
$ be the solution to Eq.(1.1) for the initial $z_{0}\in H^2\times
H^{1}$, and let $T$ be the maximal existence time of the solution
$z$, which is guaranteed by Theorem 2.2.\\

Combining (3.8), (3.11)and (3.12) we obtain
\begin{eqnarray}
&&\frac{d}{dt}\int_{\mathbb{R}}(m^2+\rho^2+\rho^2_{x})dx\\
&=&(2k_{1}-1)\int_{\mathbb{R}}m^{2}u_{x}dx+(k_{3}-k_{2})\int_{\mathbb{R}}u_{x}\rho^{2}dx\nonumber\\
&&+3k_{3}\int_{\mathbb{R}}u_{x}\rho_{x}^{2}dx
+2(k_{3}-k_{2})\int_{\mathbb{R}}u_{xx}\rho\rho_{x}dx\nonumber.
\end{eqnarray}\\

If $k_{1}\leq\frac{1}{2},$ $k_{3}\leq\min\{ k_{2},0\}$, assume that
there exist $M_{1}>0$ such that
$$u_{x}(t,x)\geq - M_{1},\ \ \ \forall \ (t,x)\in
[0,T)\times\mathbb{R}.$$  By Lemma 3.4, we have
$$\|\rho(t,\cdot)\|_{L^{\infty}}\leq
e^{-k_{3}M_{1}T}\|\rho_{0}(\cdot)\|_{L^{\infty}}, \ \ \ \forall t\in
[0,T).$$ It then follows from (3.16) that
\begin{eqnarray}
&&\frac{d}{dt}\int_{\mathbb{R}}(m^2+\rho^2+\rho^2_{x})dx\\&\leq
&(-2k_{1}+k_{2}-4k_{3}+1)M_{1}\int_{\mathbb{R}}(m^2+\rho^2+\rho_x^2)
dx\nonumber\\&&+2(k_{2}-k_{3})\|\rho(t,\cdot)\|_{L^{\infty}}\int_{\mathbb{R}}(u_{xx}^{2}+\rho^2_{x})
dx\nonumber\\&\leq
&c\int_{\mathbb{R}}(m^2+\rho^2+\rho^2_{x})dx\nonumber,
\end{eqnarray}
where
$$c=(-2k_{1}+k_{2}-4k_{3}+1)M_{1}+2(k_{2}-k_{3})e^{-k_{3}M_{1}T}\|\rho_{0}(\cdot)\|_{L^{\infty}}.$$
By means of Gronwall's inequality, we obtain $\forall \ t\in [0,T)$
\begin{eqnarray*}
\|u(t,\cdot)\|^2_{H^2}+\|\rho(t,\cdot)\|^2_{H^1}&\leq&\|m(t,\cdot)\|^2_{L^{2}}
+\|\rho(t,\cdot)\|^2_{H^1}\\
&\leq&(\|m(0,\cdot)\|^2_{L^{2}}+\|\rho(0,\cdot)\|^2_{H^1 })e^{ct}.
\end{eqnarray*}
The above inequality ensures that the solution $z$ does not blow up
in finite time.\\

If $k_{1}\geq\frac{1}{2},$ $k_{3}\geq\max\{k_{2},0\}$, assume that
there exist $M_{1}>0$ such that
$$u_{x}(t,x)\leq M_{1},\ \ \ \forall \ (t,x)\in
[0,T)\times\mathbb{R}.$$  By Lemma 3.4, we have
$$\|\rho(t,\cdot)\|_{L^{\infty}}\leq
e^{k_{3}M_{1}T}\|\rho_{0}(\cdot)\|_{L^{\infty}}, \ \ \ \forall t\in
[0,T).$$ The left proof is similar to the proof above, so we omit it
here.\\

Otherwise, let $T< \infty.$ Assume that there exists $M_{1}>0$  such
that
$$
\|u_{x}(t,\cdot)\|_{L^{\infty}}\leq M_{1},\ \ \ \forall \ t\in
[0,T).$$
 By Lemma 3.4, we have
$$\|\rho(t,\cdot)\|_{L^{\infty}}\leq
e^{|k_{3}|M_{1}T}\|\rho_{0}(\cdot)\|_{L^{\infty}}, \ \ \ \forall \
t\in [0,T).$$ It then follows from (3.16) that
\begin{eqnarray}
&&\frac{d}{dt}\int_{\mathbb{R}}(m^2+\rho^2+\rho^2_{x})dx\\&\leq
&(|2k_{1}-1|+|k_{3}-k_{2}|+3|k_{3}|)M_{1}\int_{\mathbb{R}}(m^2+\rho^2+\rho_x^2)
dx\nonumber\\&&+2|k_{2}-k_{3}|\|\rho(t,\cdot)\|_{L^{\infty}(\mathbb{R})}\int_{\mathbb{R}}(u_{xx}^{2}+\rho^2_{x})
dx\nonumber\\&\leq&(|2k_{1}-1|+|k_{3}-k_{2}|+3|k_{3}|)M_{1}\int_{\mathbb{R}}(m^2+\rho^2+\rho_x^2)
dx\nonumber\\&&+2|k_{2}-k_{3}|e^{|k_{3}|M_{1}T}\|\rho_{0}(\cdot)\|_{L^{\infty}(\mathbb{R})}\int_{\mathbb{R}}(m^{2}+\rho^2_{x})
dx\nonumber\\&\leq&c\int_{\mathbb{R}}(m^2+\rho^2+\rho^2_{x})dx\nonumber,
\end{eqnarray}
where
$c=(|2k_{1}-1|+|k_{3}-k_{2}|+3|k_{3}|)M_{1}+2|k_{2}-k_{3}|e^{|k_{3}|M_{1}T}\|\rho_{0}(\cdot)\|_{L^{\infty}}.$\\

By means of Gronwall's inequality, we obtain $\forall \ t\in [0,T)$
\begin{eqnarray*}
\|u(t,\cdot)\|^2_{H^2}+\|\rho(t,\cdot)\|^2_{H^1}&\leq&\|m(t,\cdot)\|^2_{L^{2}}
+\|\rho(t,\cdot)\|^2_{H^1}\\
&\leq&(\|m(0,\cdot)\|^2_{L^{2}}+\|\rho(0,\cdot)\|^2_{H^1})e^{ct}.
\end{eqnarray*}
The above inequality ensures that the solution $z$ does not blow up
in finite time.\\

On the other hand, by  Sobolev's imbedding theorem, we see that if
one of the conditions in the theorem holds, then the solution will
blow up in finite time. This completes the proof of
the theorem.\\
\newline
\textbf{Remark 3.1.}
Note that Theorem 3.2 shows that$$T(\|z_{0}\|_{H^{s}\times
H^{s-1}})=T(\|z_{0}\|_{H^{s^{\prime}}\times H^{s^{\prime}-1}}), \ \
\forall \ s, \ s^{\prime}>\frac{5}{2},$$ while Theorem 3.3 implies
that
$$T(\|z_{0}\|_{H^{s}\times H^{s-1}})\leq T(\|z_{0}\|_{H^{2}\times
H^{1}}),  \ \ \forall s\geq2.$$

\section{Blow up}
\noindent

In this section, we discuss the blow-up phenomena of the system
(1.1)  and prove that there exist strong solutions to (1.1) which do
not exist globally in time.\\
\newline
\textbf{Theorem 4.1.}
Let $z_{0}=\left(
                                                     \begin{array}{c}
                                                       u_{0} \\
                                                       \rho_{0} \\
                                                     \end{array}
                                                   \right)
\in H^s\times H^{s-1}, s>\frac{5}{2},$  and T be the maximal time of
the solution $z=\left(
                                    \begin{array}{c}
                                      u \\
                                      \rho\\
                                    \end{array}
                                  \right)$
to (1.1) with the initial data $z_0$. If $1<k_{1}\leq3,$ $k_{2}\geq
0,$ $u_{0}$ is odd, $\rho_{0}$ is even, $\rho_{0}(0)=0,$
$u_{0}^{\prime}(0)>0,$ then $T$ is bounded above by
$\frac{2}{k_{1}-1}\frac{1}{u_{0}^{\prime}(0)}$ and $u_{x}(t,0)$
tends to positive infinite as $t$ goes to $T.$\\
\newline
\textbf{Proof}\quad
Let $z=\left(
                                     \begin{array}{c}
                                       u \\
                                       \rho \\
                                     \end{array}
                                   \right)
$ be the solution to Eq.(1.1) for the initial $z_{0}\in H^s\times
H^{s-1}, s>\frac{5}{2}$, and let $T$ be the maximal existence time
of the solution $z$, which is guaranteed by Theorem 2.2.\\

Note that $\partial_{x}^{2}p*f=p*f-f$. Differentiating the first
equation in (2.2) with respect to $x$, then we have
\begin{equation}
u_{tx}-uu_{xx}=-\frac{k_{1}}{2}u^{2}+\frac{k_{1}-1}{2}u_{x}^{2}-\frac{k_{2}}{2}\rho^{2}
+p*(\frac{k_{1}}{2}u^{2}+\frac{3-k_{1}}{2}u_{x}^{2}+\frac{k_{2}}{2}\rho^{2}).
\end{equation}
Note that Eq.(1.1) is invariant under the transformation
$(u,x)\rightarrow(-u,-x)$ and $(\rho, x)\rightarrow(\rho, -x).$ Thus
we deduce that if $u_{0}(x)$ is odd and $\rho_{0}(x)$ is even, then
$u(t,x)$ is odd and $\rho(t,x)$ is even with respect to $x$ for any
$t\in[0,T).$ By continuity with respect to $x$ of $u$ and $u_{xx},$
we have
\begin{equation}
u(t,0)=u_{xx}(t,0)=0, \ \ \ \forall \ t\in[0,T).
\end{equation}
From Eq.(3.4) and $u$ being odd with respect to $x$, we infer that
$q(t,x)$ is odd with respect to $x$. Then we have that $q(t,0)=0$
for all $t\in[0,T).$ Since $\rho_{0}(0)=0$, it follows from Lemmas
3.3-3.4 that $\rho(t,0)=\rho(t,-k_{3}q(t,0))=0$ for all $t\in[0,T).$
Hence, in view of (4.1) and (4.2), we obtain
\begin{equation}
u_{tx}(t,0)=\frac{k_{1}-1}{2}u_{x}^{2}(t,0)
+p*(\frac{k_{1}}{2}u^{2}+\frac{3-k_{1}}{2}u_{x}^{2}+\frac{k_{2}}{2}\rho^{2})(t,0).
\end{equation}\\

By
$p*(\frac{k_{1}}{2}u^{2}+\frac{3-k_{1}}{2}u_{x}^{2}+\frac{k_{2}}{2}\rho^{2})\geq0$
and (4.3), we get
$$u_{tx}(t,0)\geq\frac{k_{1}-1}{2}u_{x}^{2}(t,0), \ \ \ t\in[0,T).$$
Set $h(t)=u_{x}(t,0).$ Since $h(0)>0,$ in view of
$u_{0}^{\prime}(0)>0,$ it follows that
\begin{equation}
0<\frac{1}{h(t)}\leq\frac{1}{h(0)}-\frac{k_{1}-1}{2}t.
\end{equation}\\
The above inequality implies that
$T<\frac{2}{k_{1}-1}\frac{1}{h(0)}$ and $u_{x}(t,0)$ tends to
positive infinite as $t$ goes to $T.$ This completes the proof of
the theorem.\\
\newline
\textbf{Theorem 4.2.}
Let $z_{0}=\left(
                                                     \begin{array}{c}
                                                       u_{0} \\
                                                       \rho_{0} \\
                                                     \end{array}
                                                   \right)
\in H^s\times H^{s-1}, s>\frac{5}{2},$  and T be the maximal time of
the solution $z=\left(
                                    \begin{array}{c}
                                      u \\
                                      \rho\\
                                    \end{array}
                                  \right)$
to (1.1) with the initial data $z_0$. If $1< k_{1}\leq3,$
 $k_{2}\geq 0,$ $u_{0}$ is odd, $\rho_{0}$ is even,
$\rho_{0}(0)=0,$ $u_{0}^{\prime}(0)=0,$ then $T$ is finite.\\
\newline
\textbf{Proof}\quad
Let $z=\left(
                                     \begin{array}{c}
                                       u \\
                                       \rho \\
                                     \end{array}
                                   \right)
$ be the solution to Eq.(1.1) for the initial $z_{0}\in H^s\times
H^{s-1}, s>\frac{5}{2}$, and let $T$ be the maximal existence time
of the solution $z$, which is guaranteed by Theorem 2.2.\\

Following the similar proof in Theorem 4.1, we have
\begin{equation}
u_{tx}(t,0)=\frac{k_{1}-1}{2}u_{x}^{2}(t,0)
+p*(\frac{k_{1}}{2}u^{2}+\frac{3-k_{1}}{2}u_{x}^{2}+\frac{k_{2}}{2}\rho^{2})(t,0).
\end{equation}
By $u_{0}^{\prime}(0)=0,$ the continuity of the ordinary
differential equation and the uniqueness, we have
$$\frac{d}{dt}u_{x}(t,0)\geq p*(\frac{k_{1}}{2}u^{2}+\frac{3-k_{1}}{2}u_{x}^{2}+\frac{k_{2}}{2}\rho^{2})(t,0)>0.$$
Therefore, $h(t)$ is strictly increasing on $[0,T).$ Since $h(0)=0,$
it follows that $h(t_{0})>0$ for some $t_{0}\in (0,T).$ Solving the
following inequality $$\frac{d}{dt}h(t)>\frac{k_{1}-1}{2}h(t)^{2},$$
we obtain
$$0<\frac{1}{h(t)}\leq\frac{1}{h(t_{0})}-\frac{k_{1}-1}{2}(t-t_{0}), \ \ \ \ t\in[t_{0}, T).$$
Consequently, we get $T<t_{0}+\frac{2}{k_{1}-1}\frac{1}{h(t_{0})}.$
This completes the proof of the theorem.\\
\newline
\textbf{Remark 4.1.}
Note that Eq.(1.1) is also invariant under the transformation
$(u,x)\rightarrow(-u,-x)$ and $(\rho, x)\rightarrow(-\rho, -x).$
Thus if the condition "$\rho_{0}$ is even, $\rho_{0}(0)=0,$" in
Theorems 4.1-4.2 and the following Corollary 4.1 is replaced by
"$\rho_{0}$ is odd" the conclusions also hold true.\\

 Then we give a corollary about Theorems 4.1-4.2.\\
\newline
\textbf{Corollary 4.1.}
Let $z_{0}=\left(
                                                     \begin{array}{c}
                                                       u_{0} \\
                                                       \rho_{0} \\
                                                     \end{array}
                                                   \right)
\in H^s\times H^{s-1}, s>\frac{5}{2},$  and T be the maximal time of
the solution $z=\left(
                                    \begin{array}{c}
                                      u \\
                                      \rho\\
                                    \end{array}
                                  \right)$
to (1.1) with the initial data $z_0$. If $1<k_{1}\leq3,$ $k_{2}\geq
0,$ $m_{0}$ is odd, $\rho_{0}$ is even, $\rho_{0}(0)=0,$
$\int_{0}^{+\infty}e^{-y}m_{0}(y)dy\geq0,$ then $T$ is finite.\\
\newline
\textbf{Proof}\quad
Note that $p(-x)=p(x),$ if $m_{0}$ is odd, then
\begin{align*}
u_{0}(x)&= \int_{\mathbb{R}}p(x-y)m_{0}(y)dy=\int_{\mathbb{R}}p(-x+y)(-m_{0}(-y))dy\\
&= \ -\int_{\mathbb{R}}p(-x-y)m_{0}(y)dy=-u_{0}(-x),
\end{align*}
from which we know that $u_{0}(x)$ is odd as well. Since
$$u_{0}^{\prime}(x)=-\frac{1}{2}e^{-x}\int_{-\infty}^{x}e^{y}m_{0}(y)dy
+\frac{1}{2}e^{x}\int_{x}^{+\infty}e^{-y}m_{0}(y)dy,$$ we get
$$u_{0}^{\prime}(0)=\int_{0}^{+\infty}e^{-y}m_{0}(y)dy.$$ It follows
from Theorems 4.1-4.2 that the corresponding solution to Eq.(1.1)
blows up infinite time.

\bigskip
\noindent\textbf{Acknowledgments} This work was partially
supported by NNSFC (No. 10971235), RFDP (No. 200805580014),
NCET-08-0579 and the key project of Sun Yat-sen University.

\end{document}